\documentclass[10pt,letterpaper]{article}
\usepackage{blindtext,titlefoot}
\usepackage{authblk}
\usepackage{marvosym}
\usepackage[utf8]{inputenc}
\usepackage[english]{babel}
\usepackage{amsmath}
\usepackage{amsfonts}
\usepackage{amssymb}
\usepackage{makeidx}
\usepackage{graphicx}
\usepackage{lmodern}
\usepackage{kpfonts}
\usepackage{dsfont}
\usepackage{cancel}
\usepackage{amsthm} 
\usepackage{hyperref}
\usepackage[left=2cm,right=2cm,top=2cm,bottom=2cm]{geometry}
\usepackage{subcaption}
\usepackage{multirow}
\usepackage{float}
\usepackage{url}
\usepackage{multicol}
\usepackage{mathtools, cuted}
\usepackage{lipsum}
\usepackage{booktabs}
\usepackage{comment}
\usepackage{cite}

\usepackage{fancyhdr}
\fancypagestyle{firstpage}{\fancyhf{}
\setcounter{page}{1}
\rfoot{\thepage}
}

\providecommand{\nset}[1]{
\mathbb{#1}
}
\providecommand{\set}[1]{
\left\{#1\right\}
}
\providecommand{\com}[1]{``#1"}

\providecommand{\ifr}[5]{
{}^{#1}_{#2}{#3}_{#4}^{#5}
}
\providecommand{\gam}[1]{
\Gamma\left(#1 \right)
}

\providecommand{\norm}[1]{
\left\lVert #1 \right\rVert
}

\providecommand{\ds}[1]{
\displaystyle #1
}
\providecommand{\der}[3]{
\dfrac{#1^{#3} }{ #1 #2^{#3}}
}
\providecommand{\dt}[1]{
\mbox{det}\left(#1\right)
}

\newtheorem{theorem}{ Theorem}[section]
\newtheorem{definition}[theorem]{Definition}
\newtheorem{proposition}[theorem]{Proposition}
\newtheorem{example}[theorem]{Example}

\usepackage{enumitem} 
\setlist[itemize]{noitemsep} 

\usepackage{titlesec} 
\titleformat{\section}[block]{\large\bfseries\scshape\centering}{\thesection.}{1em}{} 
\titleformat{\subsection}[block]{\large\bfseries\scshape\centering}{\thesubsection.}{1em}{}
\titleformat{\subsubsection}[block]{\large\bfseries\scshape\centering}{\thesubsubsection.}{1em}{} 

\title{\huge\bfseries Fractional pseudo-Newton method and its use in the solution of a nonlinear system  that allows the construction of a hybrid solar receiver}

\author[,a]{\normalsize A. Torres-Hernandez  \footnote{Email address: anthony.torres@ciencias.unam.mx; Corresponding author; ORCID: 0000-0001-6496-9505}}
\affil[a]{\normalsize Department of Physics, Faculty of Science - UNAM, Mexico}

\author[,b]{\normalsize F. Brambila-Paz \footnote{Email address: fernandobrambila@gmail.com; ORCID: 0000-0001-7896-6460}}
\affil[b]{\normalsize Department of Mathematics, Faculty of Science - UNAM, Mexico}

\affil[c]{\normalsize Faculty of Engineering, Universidad Panamericana - Aguascalientes, Mexico}

\author[,c,d]{\normalsize P. M. Rodrigo \footnote{Email address: prodrigo@up.edu.mx; ORCID: 0000-0003-0100-6124}}
\affil[d]{\normalsize Centre for Advanced Studies on Energy and Environment (CEAEMA), University of Jaén, Spain.}

\author[,c]{\normalsize E. De-la-Vega \footnote{Email address: evega@up.edu.mx; ORCID: 0000-0001-9491-6957}}

\author[,c]{\normalsize C. C. Calabrese \footnote{Email address: ccalabrese@up.edu.mx; ORCID: 0000-0001-9844-3368}}

\date{}

\begin{document}

\maketitle
\thispagestyle{firstpage}


\begin{abstract}

The following document presents a possible solution and a brief stability analysis for a nonlinear system, which is obtained by studying the possibility of building a hybrid solar receiver; It is necessary to mention that the solution of the aforementioned system is relatively difficult to obtain through iterative methods since the system is apparently unstable.  To find this possible solution is used a novel numerical method valid for one and several variables, which using the fractional derivative, allows us to find solutions for some nonlinear systems in the complex space using real initial conditions, this method is also valid for linear systems. The method described above has an order of convergence (at least) linear, but it is easy to implement and it is not necessary to invert some matrix for solving nonlinear systems and linear systems.

\end{abstract}

\section{Introduction}

A classic problem of common interest in Physics, Mathematics and Engineering is to find the zeros of a function $f:\Omega \subset \nset{R}^n \to \nset{R}^n$, that is,

\begin{eqnarray*}
\set{\xi \in \Omega \ : \ \norm{f(\xi)}=0},
\end{eqnarray*}

this problem often arises as a consequence of wanting to solve other problems, for instance, if we want to determine the eigenvalues of a matrix or want to build a box with a given volume but with a minimal surface; in the first example, we need to find the zeros (or roots) of the characteristic polynomial of the matrix, while in the second one we need to find the zeros of the gradient of a function that relates the surface of the box with its volume.

Although finding the zeros of a function may seem like a simple problem, in general, it involves solving nonlinear equations and numerical methods are needed to try to determine the solutions to these problems; it should be noted that when using numerical methods, the word \com{determine} should be interpreted as to approach a solution with a degree of precision desired. The numerical methods mentioned above are usually of the iterative type and work as follows: suppose we have a function $f : \Omega \subset \nset{R}^n \to \nset{R}^n$ and we search a value $\xi\in \nset{R}^n$ such that $\norm{f (\xi)} = 0$, then we may start by giving an initial value $x_0\in \nset{R}^n$ and then calculate a value $x_i$ close to the searched value $\xi$ using an iteration function $\Phi : \nset{R}^n \to \nset{R}^n$ as follows 

\begin{eqnarray}\label{eq:c2.01}
x_{i+1}:=\Phi(x_i), & i=0,1,2,\cdots.
\end{eqnarray} 

When it is assumed that the iteration function $ \Phi $ is continuous around $ \xi $ and that the sequence $ \set{x_i} _{i = 0} ^ \infty $ converges to $ \xi $, it holds that

\begin{eqnarray}\label{eq:c2.07}
\xi=\lim_{i\to \infty}x_{i+1}=\lim_{i\to \infty}\Phi(x_i)=\Phi\left(\lim_{i\to \infty}x_i\right)=\Phi(\xi),
\end{eqnarray}

the previous result is the reason why the method given in \eqref{eq:c2.01} is called \textbf{fixed point method}.

In the last section of this document, we study the nonlinear system that describes a hybrid solar panel, which consists of a photovoltaic-thermoelectric generator, and we will proceed to find a possible solution for this system using the fractional pseudo-Newton method, because the apparent instability of the system makes the classic Newton's method not the most suitable to solve it.

\section{Previous works}

\subsection{Historical background of fractional calculus}

The question that led to the emergence of a new branch of mathematical analysis known as fractional calculus, was asked by L'H\^opital in 1695 in a letter to Leibniz, as a consequence of the notation $ d ^ nf (x) / dx ^ n $; perhaps it was a game of symbols that which prompted L'H\^opital to ask Leibniz: \com {What happens if $ n = 1/2 $?}, Leibniz replied in a letter, almost prophetically: \com{$ \cdots $ is an apparent paradox from which, one day, useful consequences will be drawn \cite{miller93}.}  Subsequently, the question became: may the order $ n $ of the derivative  be any number: rational, irrational or complex? Because the question was answered affirmatively, the name of the fractional calculus has become an incorrect name and it would be more correct to call it \textbf{arbitrary order integration and differentiation}.

The concepts of arbitrary order differentiation and integration are not new. Interest in these subjects was evident almost in tandem with the emergence of conventional calculus (differentiation and integration of the integer order), the first systematic studies were written in the early and mid-19th century by Liouville (1832), Riemann (1953), and Holmgrem (1864), although Euler (1730), Lagrange (1772), and other authors made contributions even earlier \cite{miller93}.

When a function does not have integer derivative the notion of weak derivative is required. The weak derivatives give rise to generalized functions or distributions, which are often used in quantum mechanics. It is important to mention that there are functions that do not have weak derivative but have fractional derivative, such as the Cantor function \cite{nigmatullin1992fractional}.

Caputo in 1967, developed the first application of fractional calculus related to diffusion processes, in what he named as anomalous diffusion equation \cite{kilbas2006theory}. There is practically no branch of classical analysis that remains exempt from fractional calculus.

\subsection{Fractional Iterative Methods}

Although the interest in fractional calculus has mainly focused on the study and development of techniques to solve systems of differential equations of non-integer order \cite{hilfer00,kilbas2006theory, martinez2017applications1,martinez2017applications2,torreshern2019proposal}, over the years, iterative methods have also been developed that use the properties of fractional derivatives to solve algebraic equation systems \cite{gao2009local,fernando2017fractional,brambila2018fractional,akgul2019fractional,torreshern2020}.
These methods in general may be called \textbf{fractional iterative methods}.

It should be mentioned that depending on the nature of the definition of fractional derivative used, fractional iterative methods have the particularity that they may be used of local form \cite{gao2009local} or of global form \cite{torreshern2020}. These methods also allow searching for complex roots for polynomials using only real initial conditions.

Although using an iterative method that uses fractional derivatives seems to require an unnecessary effort, considering that a more natural option would be Newton's method \cite{ortega1970iterative}, it should be noted that Newton's method may only be used in differentiable functions, while fractional derivatives may be used in a larger number of functions \cite{nigmatullin1992fractional}. Some differences between Newton's method and two fractional methods are listed in the Table \ref{tab:06}

\begin{table}[!ht]
\centering
\footnotesize
\begin{tabular}{c|c|c|c}
\toprule
&Classic Newton&Fractional Newton& Fractional Pseudo-Newton \\ \midrule
\begin{tabular}{c}
Can it find the complex zeros\\
 of a polynomial using \\
 real initial conditions?
\end{tabular}& No & Yes & Yes \\ \midrule
\begin{tabular}{c}
Can it find multiple zeros \\
of a function using a \\
single initial condition?
\end{tabular}& No & Yes & Yes \\ \midrule
\begin{tabular}{c}
It can be used if the function \\ 
is not differentiable?
\end{tabular}& No&Yes&Yes \\ \midrule
\begin{tabular}{c}
For a space of dimension $ N $ \\
are needed
\end{tabular}&
\begin{tabular}{c}
$ N \times N $ classic \\
 partial derivatives
\end{tabular}& 
\begin{tabular}{c}
$ N \times N $ fractional \\
 partial derivatives
\end{tabular}
 & \begin{tabular}{c}
$ N$ fractional \\
 partial derivatives
\end{tabular} \\ \midrule
\begin{tabular}{c}
Is it recommended to solve systems\\
 where the partial derivatives are \\
 analytically difficult to obtain?
\end{tabular}& Not recommended& Not recommended &Is recommended \\ \bottomrule
\end{tabular}
\caption{Some differences between the classical Newton method and two fractional iterative methods.}\label{tab:06}
\end{table}

\subsection{Introduction of a Hybrid Solar Receiver}

Concentrator photovoltaic (CPV) systems represent a technological success in solar energy applications because of the high conversion efficiencies commercially achieved. These systems use optical devices to concentrate the sunlight onto small highly efficient multi-junction (MJ) solar cells. Efficiency records in a laboratory of $47.1\%, \ 43.4\%$, and $38.9\%$  at cell, mono-module, and module level respectively have been shown \cite{philipps2015current,green2019solar}, while commercial CPV modules have a mean efficiency of $30.0\%$ \cite{perez2018efficiency}. In spite of these high efficiencies, CPV systems have a higher levelized cost of energy (LCOE) than traditional photovoltaic (PV) systems \cite{kost2013levelized,talavera2016worldwide}. Among the strategies that are being investigated to make CPV systems more competitive, the increase of the concentration factor or the increase of efficiency are considered.

An increase of efficiency can be achieved by recovering part of the waste heat generated in the solar cells. Among the strategies to get this, the hybridization with thermoelectric generators (TEG) is being proposed. Thermoelectric (TE) materials such as Bi${}_2$Te${}_3$ can operate under the Seebeck effect to transform a heat flux to electricity \cite{rowe2018crc}. Hybrid CPV-TEG systems are in the research stage and several laboratory prototypes have been reported \cite{beeri2015hybrid,tamaki2017hybrid,kil2017highly}, although currently they are far from the expected benefits.

\section{Fractional Pseudo-Newton Method}

\subsection{Order of Convergence}

Before continuing it is necessary to have the following definition \cite{plato2003concise}

\begin{definition}
Let $ \Phi: \Omega \subset \nset{R}^ n \to \nset{R}^ n $ be an iteration function with a fixed point $ \xi \in \Omega $. Then the method \eqref{eq:c2.01} is called (locally) \textbf{convergent of (at least) order $ \boldsymbol{p} $} ($ p \geq 1 $), if exists $ \delta> 0 $  and exists a non-negative constant $ C $ (with $ C <1 $ if $ p = 1 $) such that for any initial value $ x_0 \in B (\xi; \delta) $ it holds that

\begin{eqnarray}\label{eq:c2.08}
\norm{x_{k+1}-\xi}\leq C \norm{x_k-\xi}^p, & k=0,1,2,\cdots,
\end{eqnarray}

where $ C $ is called convergence factor.

\end{definition}

The order of convergence is usually related to the speed at which the sequence generated by \eqref{eq:c2.01} converges. For the particular cases $ p = 1 $ or $ p = 2 $ it is said that the method has (at least) linear convergence or (at least) quadratic convergence, respectively. The following theorem \cite{plato2003concise,stoer2013}, allows characterizing the order of convergence of an iteration function $ \Phi $ with its derivatives 

\begin{theorem}\label{teo:c2.01}
Let $ \Phi: \Omega \subset \nset{R}^n \to \nset {R}^n $ be an iteration function with a fixed point $ \xi \in \Omega $. Assuming that $\Phi $ is $ p$-times differentiable in $ \xi $ for some $ p \in \nset{N} $, and furthermore

\begin{eqnarray}\label{eq:c2.09}
\left\{
\begin{array}{cc}
 \norm{\Phi^{(k)}(\xi)}=0, \ \forall  k\leq p-1, & \mbox{if }p\geq 2 \\
 \norm{\Phi^{(1)}(\xi)}<1, & \mbox{if }p=1
\end{array}\right.,
\end{eqnarray}

then $ \Phi $ is (locally) convergent of (at least) order $ p $.

\end{theorem}

The previous theorem is usually very useful to generate a fixed point method with an order of convergence desired, an order of convergence that is usually appreciated in iterative methods is the (at least) quadratic order. If we have a function $f:\Omega \subset \nset{R}^n \to \nset{R}^n$ and we search a value $\xi\in \Omega$ such that $\norm{f(\xi)}=0$, we may build an iteration function $ \Phi $ in general form as \cite{burden2002analisis}

\begin{eqnarray}\label{eq:c2.13}
\Phi(x)=x-A(x)f(x),
\end{eqnarray}

with $
A(x):=\left([A]_{jk}(x) \right)$ a matrix, where $[A]_{jk}:\nset{R}^n \to \nset{R}$  ($1\leq j,k\leq n$). Notice that the matrix $ A (x) $ is determined according to the order of convergence desired. Denoting by $ \dt{A} $ the determinant of the  matrix $ A $, is possible to demonstrate that any matrix $ A (x) $ that fulfill the following condition \cite{torreshern2020}

\begin{eqnarray}\label{eq:c2.17}
\lim_{x\to \xi}A(x)= \left(f^{(1)}(\xi)\right)^{-1}, & \dt{f^{(1)}(\xi)}\neq 0,
\end{eqnarray}

where $f^{(1)}$ is the \textbf{Jacobian matrix} of the function $f$ \cite{ortega1990numerical}, guarantees that $\norm{\Phi^{(1)}(\xi)}=0$. As a consequence, exists $ \delta> 0 $ such that the iteration function $ \Phi $ given by \eqref{eq:c2.13}, converges (locally) with an order of convergence (at least) quadratic in $ B (\xi; \delta) $.

\subsection{Fractional Derivative of Riemann-Liouville}

One of the key pieces in the study of fractional calculus is the iterated integral, which is defined as follows \cite{hilfer00}

\begin{definition}
Let $ L_{loc} ^ 1 (a, b) $, the space of locally integrable functions in the interval $ (a, b) $. If $ f $ is a function such that $ f \in L_ {loc} ^ 1 (a, \infty) $, then the $n$-th iterated integral of the function $ f $ is given by \cite{hilfer00}

\begin{eqnarray}\label{eq:c1.16}
\begin{array}{c}
\ds \ifr{}{a}{I}{x}{n} f(x)=\ifr{}{a}{I}{x}{}\left(\ifr{}{a}{I}{x}{n-1} f(x)  \right)=\frac{1}{(n-1)!}\int_a^x(x-t)^{n-1}f(t)dt,
\end{array}
\end{eqnarray}

where

\begin{eqnarray*}
\ifr{}{a}{I}{x}{} f(x):=\int_a^x f(t)dt.
\end{eqnarray*}

\end{definition}

Considerate that $ (n-1)! = \gam{n} $
, a generalization of \eqref{eq:c1.16} may be obtained for an arbitrary order $ \alpha> 0 $

\begin{eqnarray}\label{eq:c1.17}
\ifr{}{a}{I}{x}{\alpha} f(x)=\dfrac{1}{\gam{\alpha}}\int_a^x(x-t)^{\alpha-1}f(t)dt,
\end{eqnarray}

the equations \eqref{eq:c1.17}  correspond to the definitions of \textbf{(right) fractional integral of Riemann-Liouville}.  Fractional integrals satisfy the  \textbf{semigroup property}, which is given in the following proposition \cite{hilfer00}

\begin{proposition}
Let $ f $ be a function. If $ f \in L_{loc} ^ 1 (a, \infty) $, then the fractional integrals of $ f $ satisfy that

\begin{eqnarray}\label{eq:c1.19}
\ifr{}{a}{I}{x}{\alpha} \ifr{}{a}{I}{x}{\beta}f(x) = \ifr{}{a}{I}{x}{\alpha + \beta}f(x),& \alpha,\beta>0.
\end{eqnarray}

\end{proposition}

From the previous result, and considering that the operator $ d / dx $  is the inverse operator to the left of the operator $ \ifr {}{a}{I}{x}{} $, any integral $ \alpha$-th of a function $ f \in L_{loc} ^ 1 (a, \infty) $ may be written as

\begin{eqnarray}\label{eq:c1.20}
\ifr{}{a}{I}{x}{\alpha}f(x)=\dfrac{d^n}{dx^n}\ifr{}{a}{I}{x}{n}\left( \ifr{}{a}{I}{x}{\alpha}f(x) \right)=\dfrac{d^n}{dx^n}\left( \ifr{}{a}{I}{x}{n+\alpha}f(x)\right).
\end{eqnarray}

With the previous results, we can built the operator  \textbf{fractional derivative of Riemann-Liouville}, as follows \cite{hilfer00,kilbas2006theory}

\begin{eqnarray}\label{eq:c1.23}
\normalsize
\begin{array}{c}
\ifr{}{a}{D}{x}{\alpha}f(x) := \left\{
\begin{array}{cc}
\ds \ifr{}{a}{I}{x}{-\alpha}f(x), &\mbox{if }\alpha<0\\  
\ds \dfrac{d^n}{dx^n}\left( \ifr{}{a}{I}{x}{n-\alpha}f(x)\right), & \mbox{if }\alpha\geq 0
\end{array}
\right.
\end{array}, 
\end{eqnarray}

where  $ n = \lfloor \alpha \rfloor + 1 $. Applying the  operator \eqref{eq:c1.23} with $ a = 0 $ and $ \alpha \in \nset{R}\setminus\nset {Z} $ to the  function $ x^{\mu} $, with $\mu>-1$, we obtain that

\begin{eqnarray}\label{eq:c1.13}
\ifr{}{0}{D}{x}{\alpha}x^\mu = 
 \dfrac{\gam{\mu+1}}{\gam{\mu-\alpha+1}}x^{\mu-\alpha}.
\end{eqnarray}

\subsection{Iteration Function of Fractional Pseudo-Newton Method}

Let $f$ a function, with $f:\Omega \subset \nset{R}\to \nset{R}$. We can consider the problem of finding a value $\xi\in \nset{R}$ such that $f(\xi)$=0. A first approach to value $\xi$ is by a linear approximation of the function $f$ in a valor $x_i \approx \xi$, that is,

\begin{eqnarray}\label{eq:c2.36}
f(x)\approx f(x_i)+f^{(1)}(x_i)(x-x_i),
\end{eqnarray}

then, assuming that $ \xi $ is a zero of $ f $, from the previous expression we have that

\begin{eqnarray*}
\xi \approx x_i- \left(f^{(1)}(x_i) \right)^{-1} f(x_i),
\end{eqnarray*}

as consequence, a sequence $ \set{x_i}_{i = 0} ^ \infty $ that approximates the value $ \xi $ may be generated using the iteration function

\begin{eqnarray*}
x_{i+1}:=\Phi(x_i)=x_i- \left(f^{(1)}(x_i) \right)^{-1} f(x_i), & i=0,1,2,\cdots.
\end{eqnarray*}

which corresponds to well-known \textbf{Newton's method} \cite{ortega1970iterative}. However, the equation \eqref{eq:c2.36} is not the only way to generate a linear approximation to the  function $ f $ in the point $ x_i $, in general it may be taken as

\begin{eqnarray}\label{eq:c2.37}
f(x)\approx f(x_i)+m(x-x_i),
\end{eqnarray}

where $ m $ is any constant value of a slope, that allows the approximation  \eqref{eq:c2.37} to the  function $ f $ to be valid. The previous equation allows to obtain the following iteration function

\begin{eqnarray}\label{eq:c2.38}
x_{i+1}:=\Phi(x_i)=x_i- m^{-1}f(x_i),&  i=0,1,2,\cdots,
\end{eqnarray}

which originates the \textbf{parallel chord method} \cite{ortega1970iterative}. The iteration function \eqref{eq:c2.38} can be generalized to larger dimensions as follows

\begin{eqnarray}\label{eq:c2.39}
x_{i+1}:=\Phi(x_i)= x_i- \left( m^{-1}I_n \right) f(x_i), & i=0,1,2\cdots,
\end{eqnarray}

where $ I_n $ corresponds to the identity matrix of $n\times n$. It should be noted that the idea behind the parallel chord method in several variables is just to apply \eqref{eq:c2.38} component by component.

Before continuing it is necessary to mention that for some definitions of the fractional derivative, it is satisfied that the derivative of the order $ \alpha $ of a constant is different from zero, that is,

\begin{eqnarray}\label{eq:c2.30}
\partial_k^\alpha c :=\der{\partial}{[x]_k}{\alpha}c \neq 0 , & c=constant,
\end{eqnarray}

where $ \partial_k ^ \alpha $ denotes any fractional derivative applied only in the component $ k $, that does not cancel the constants and that satisfies the following continuity relation with respect to the order $ \alpha $ of the derivative

\begin{eqnarray}\label{eq:c2.301}
\lim_{\alpha \to 1}\partial_k^\alpha c=\partial_kc.
\end{eqnarray}

Using as a basis the idea of \eqref{eq:c2.39}, and considering any fractional derivative that satisfies the conditions \eqref{eq:c2.30} and \eqref{eq:c2.301}, we can define the \textbf{fractional pseudo-Newton method} as follows

\begin{eqnarray}\label{eq:c2.40}
x_{i+1}:=\Phi(\alpha, x_i)= x_i- P_{\epsilon,\beta}(x_i) f(x_i), & i=0,1,2\cdots,
\end{eqnarray}

with $\alpha\in[0,2]\setminus\nset{Z}$, where $ P_{\epsilon, \beta} (x_i) $ is a matrix evaluated in the value $ x_i $, which is given as follows

\begin{eqnarray}
P_{\epsilon,\beta}(x_i):=\left([P_{\epsilon,\beta}]_{jk}(x_i)\right)=\left( \partial_k^{\beta(\alpha,[x_i]_k)}\delta_{jk}+ \epsilon\delta_{jk}  \right)_{x_i},
\end{eqnarray}

where

\begin{eqnarray}
\partial_k^{\beta(\alpha,[x_i]_k)}\delta_{jk}:= \der{\partial}{[x]_k}{\beta(\alpha,[x_i]_k)}\delta_{jk}, & 1\leq j,k\leq n,
\end{eqnarray}

with $ \delta_{jk} $ the Kronecker delta, $ \epsilon $ a positive constant $ \ll 1 $, and $ \beta (\alpha, [x_i]_k) $ defined as follows

\begin{eqnarray}\label{eq:c2.34}
\beta(\alpha,[x_i]_k):=\left\{
\begin{array}{cc}
\alpha, &\mbox{if }  |[x_i]_k|\neq 0\\
1,& \mbox{if }  |[x_i]_k|=0
\end{array}\right.,
\end{eqnarray}

It should be mentioned that the value $ \alpha = 1 $ in \eqref{eq:c2.34}, is taken to avoid the discontinuity that is generated when using the fractional derivative of constants in the value $ x = 0 $. 

Since in the fractional pseudo-Newton method, the matrix $P_{\epsilon,\beta}(x_i)$  does not satisfy the condition \eqref{eq:c2.17}, any sequence $ \set{x_i} _ {i = 0} ^ \infty $ generated by the iteration function \eqref {eq:c2.40} has at most one order of convergence (at least) linear. 

\subsubsection{Some Examples}

\begin{example}

Let the function:

\begin{eqnarray*}
\begin{array}{c}
f(x)= \left(\dfrac{1}{2}\sin(x_1x_2)-\dfrac{x_2}{4\pi}-\dfrac{x_1}{2}, \left( 1-\dfrac{1}{4\pi}\right)\left(e^{2x_1}-e \right)+\dfrac{e}{\pi}x_2-2ex_1 \right)^T,
\end{array}
\end{eqnarray*}

then the value $x_0=(1.03,1.03)^T$ is chosen to use the iteration function given by \eqref{eq:c2.40}, and using the fractional derivative given by \eqref{eq:c1.13}, we obtain the results of the Table \ref{tab:03}

\begin{table}[!ht]
\centering
$
\begin{array}{c|cccccc}
\toprule
&\alpha_m& {}^m\xi_1& {}^m\xi_2 &\norm{{}^m \xi - {}^{m-1} \xi }_2  &\norm{f\left({}^m\xi \right)}_2& R_m \\ 
\midrule
1	&	0.78562	&	 1.03499277 - 0.53982128i	&	5.41860852 + 4.04164098i	&	5.62354e-6	&	8.38442e-5	&	66	\\
2	&	 0.78987	&	  0.29945564        	&	  2.83683317        	&	1.09600e-5	&	9.63537e-5	&	  88	\\
3	&	 0.82596	&	 -0.26054499        	&	  0.62286899        	&	5.66073e-5	&	9.87374e-5	&	 140	\\
4	&	 0.82671	&	  -0.1561964 - 1.02056003i	&	  2.26280132 - 5.71855964i	&	4.32875e-6	&	9.51178e-5	&	 194	\\
5	&	 0.83158	&	  1.03499697 + 0.53981525i	&	  5.41862187 - 4.04161017i	&	3.94775e-6	&	8.80344e-5	&	  84	\\
6	&	 0.85861	&	  1.16151359 - 0.69659512i	&	  8.27130854 + 6.3096935i 	&	2.14707e-6	&	9.38721e-5	&	 164	\\
7	&	 1.15911	&	  1.48131686        	&	 -8.38362876        	&	1.20669e-6	&	9.56674e-5	&	 191	\\
8	&	 1.24977	&	 -1.10844524 + 0.10906317i	&	 -4.18608959 + 0.66029327i	&	3.71508e-6	&	9.81146e-5	&	 164	\\
9	&	 1.25662	&	 -1.10844605 - 0.10906368i	&	 -4.18608629 - 0.66029181i	&	3.69483e-6	&	9.66271e-5	&	 170	\\
10	&	 1.26128	&	  1.33741853  	&	 -4.14026671  	&	1.89913e-5	&	8.51053e-5	&	  67	\\
\bottomrule
\end{array}
$
\caption{Results obtained using the iterative method \eqref{eq:c2.40} with $\epsilon=-3$.}\label{tab:03}
\end{table}

\end{example}

\begin{example}

Let the function:

\begin{eqnarray*}
\normalsize
\begin{array}{c}
f(x)= \left(  -3.6 x_3\left( x_1^3x_2 +1 \right) - 3.6\cos\left( x_2^2 \right) + 10.8, 
-1.6x_1\left( x_1 +x_2^3x_3\right) - 1.6\sinh\left( x_3\right) + 6.4, 
-4x_2\left( x_1x_3^3 +1 \right) - 4\cosh\left( x_1 \right) + 24 \right)^T,
\end{array}
\end{eqnarray*}

then the value $x_0=(1.12,1.12,1.12)^T$ is chosen to use the iteration function given by \eqref{eq:c2.40}, and using the fractional derivative given by \eqref{eq:c1.13}, we obtain the results of the Table \ref{tab:04}

\begin{table}[!ht]
\centering
\footnotesize
$
\begin{array}{c|ccccccc}
\toprule
&\alpha_m& {}^m\xi_1& {}^m\xi_2 &{}^m\xi_3&\norm{{}^m \xi - {}^{m-1} \xi }_2  &\norm{f\left({}^m\xi \right)}_2& R_m \\ 
\midrule
1	&	0.96743	&	0.38147704 + 1.10471108i	&	-0.43686196 - 1.3473184i 	&	-0.38512615 - 1.4903386i 	&	2.41936e-6	&	7.07364e-5	&	57	\\
2	&	 0.96745	&	 -0.78311553 + 0.96791081i	&	 -0.58263802 + 1.2592471i 	&	  0.18175185 - 1.49135484i	&	3.85644e-6	&	8.58385e-5	&	 37	\\
3	&	 0.96766	&	  0.71500126 - 1.02632085i	&	  0.53575431 - 1.314774i  	&	  0.45273307 - 1.35710557i	&	3.01511e-6	&	8.34643e-5	&	 41	\\
4	&	 0.9677 	&	 -0.34118928 + 1.19432023i	&	  0.37199268 - 1.40125985i	&	 -0.63215137 + 1.3074313i 	&	2.72698e-6	&	8.69377e-5	&	 49	\\
5	&	 0.96796	&	  0.71500155 + 1.0263218i 	&	  0.53575489 + 1.31477495i	&	  0.45273303 + 1.35710453i	&	2.52069e-6	&	7.11216e-5	&	 34	\\
6	&	 0.97142	&	 -0.34118945 - 1.19432007i	&	  0.37199303 + 1.40125973i	&	 -0.63215109 - 1.3074314i 	&	2.32465e-6	&	8.66652e-5	&	 61	\\
7	&	 0.9718 	&	  0.38147878 - 1.10471296i	&	 -0.43686073 + 1.34732029i	&	  -0.3851262 + 1.49033466i	&	2.38466e-6	&	8.53472e-5	&	 50	\\
8	&	 0.97365	&	 -0.78311508 - 0.96791138i	&	 -0.58263753 - 1.2592475i 	&	  0.18175161 + 1.49135503i	&	1.99078e-6	&	7.57517e-5	&	 51	\\
9	&	 1.03148	&	  1.34508926      	&	 -1.29220278      	&	 -1.44485467      	&	3.68616e-6	&	9.58451e-5	&	 59	\\
10	&	 1.04155	&	 -1.43241693 	&	  1.27535274  	&	 -1.11183615 	&	4.06891e-6	&	8.95830e-5	&	 48	\\
\bottomrule
\end{array}
$
\caption{Results obtained using the iterative method \eqref{eq:c2.40} with $\epsilon=e-3$.}\label{tab:04}
\end{table}

\end{example}

When we work with a linear system of the form

\begin{eqnarray*}
Ax=b,
\end{eqnarray*}

it is possible to solve it using the method \eqref{eq:c2.40} considering the function

\begin{eqnarray*}
f(x)=Ax-b.
\end{eqnarray*}

\begin{example}

Let the function:

\begin{eqnarray*}
\begin{array}{c}
f(x)= \left(  5x_1-4x_2+3x_3-18,2x_1+5x_2-6x_3-24,-2x_1+7x_2+12x_3-30 \right)^T,
\end{array}
\end{eqnarray*}

then the value $x_0=(0.64,0.64,0.64)^T$ is chosen to use the iteration function given by \eqref{eq:c2.40}, and using the fractional derivative given by \eqref{eq:c1.13}, we obtain the results of the Table \ref{tab:05}

\begin{table}[!ht]
\centering
$
\begin{array}{c|ccccccc}
\toprule
&\alpha_m& {}^m\xi_1& {}^m\xi_2 &{}^m\xi_3 &\norm{{}^m \xi - {}^{m-1} \xi }_2  &\norm{f\left({}^m\xi \right)}_2& R_m \\ 
\midrule
 1	&	   0.90162	&	5.97144261	&	3.88571164	&	1.22857594	&	2.55098e-6	&	9.53968e-5	&	65	\\
\bottomrule
\end{array}
$
\caption{Results obtained using the iterative method \eqref{eq:c2.40} with $\epsilon=e-3$.}\label{tab:05}
\end{table}

\end{example}

\section{Equations of a Hybrid Solar Receiver}

Considering the notation

\begin{eqnarray*}
X=(x,y,z,v,w)^T:=(T_{cell},T_{hot},T_{cold},\eta_{cell},\eta_{TEG})^T,
\end{eqnarray*}

it is possible to define the following system of equations that corresponds to the combination of a solar photovoltaic system with a thermoelectric generator system \cite{bjork2015performance,bjork2018maximum}, which is named as a \textbf{hybrid solar receiver} \cite{li2018review}

\begin{eqnarray}\label{eq:004}
\left\{
\begin{array}{l}
x=y+a_1\cdot a_2(1-v)\\
y=z+a_1\cdot a_3 (1-v)(1-w)\\
z=a_4+a_1\cdot a_5 (1-v)(1-w)\\
v=a_6x+a_7\\
w=(a_8-1)\left(1-\dfrac{z+a_9}{y+a_9} \right)\left(a_8+ \dfrac{z+a_9}{y+a_9}\right)^{-1}
\end{array}\right.,
\end{eqnarray}

whose deduction and some details about the difficulty in finding its solution may be found in the reference \cite{rodrigo2019performance}. The $ a_i $'s in the previous systems are constants defined by the following expressions

\begin{eqnarray*}
\left\{
\begin{array}{l}
\begin{array}{l}
 a_2=r_{cell}+r_{sol}+A_{cell}\left(\dfrac{r_{cop}+r_{cer}}{A_{TEG}}+\dfrac{r_{intercon}}{0.5\cdot \sqrt{f^*\cdot A_{TEG}}\left(b\cdot \sqrt{f^*}+\sqrt{A_{TEG}} \right) } \right)\\
a_5=A_{cell}\left( \dfrac{r_{intercon}}{0.5\cdot \sqrt{f^*\cdot A_{TEG}}\left(b\cdot \sqrt{f^*}+\sqrt{A_{TEG}} \right) }+\dfrac{r_{cer}}{A_{TEG}}+R_{heat\_ exch} \right)
\end{array}\\
\begin{array}{lll}
a_1=\eta_{opt}\cdot C_g \cdot DNI,&
a_3=\dfrac{A_{cell}\cdot l}{f^*\cdot A_{TEG}\cdot k_{TEG}},&a_4=T_{air} \\ \\
a_6=-\eta_{cell,ref}\cdot \gamma_{cell},&a_7=\eta_{cell,ref}\left(1+25 \cdot \gamma_{cell} \right),&
a_8=\sqrt{1+ZT}\\ \\
a_9=273.15&
\end{array}
\end{array}\right.,
\end{eqnarray*}

with the following particular values \cite{rodrigo2019performance}

\begin{eqnarray*}
\left\{
\begin{array}{lll}
\eta_{opt}=0.85, & r_{intercon}=2.331 e-7 ,& T_{air}=20 \\
        C_g=800, &     A_{cell}=9e-6, &     R_{heat\_exch}=0.5 \\
DNI=900 ,&     A_{TEG}=5.04e-5, &     \eta_{cell,ref}=0.43 \\
r_{cell}=3e-6, &     f^*=0.7, &     \gamma_{cell}=4.6e-4 \\
r_{sol}=1.603e-6, & b=5e-4, &     ZT=1 \\
r_{cop}=7.5e-7, &     l=5e-4 ,& r_{cer}=8e-6\\
k_{TEG}=1.5 
\end{array}
\right..
\end{eqnarray*}

Using the system of equations \eqref{eq:004}, it is possible to define a  function $f:\Omega \subset \nset{R}^5\to \nset{R}^5$, that is,

\begin{eqnarray}\label{eq:005}
f(X):=\begin{pmatrix}
-x+y+a_1\cdot a_2(1-v)\\
-y+z+a_1\cdot a_3 (1-v)(1-w)\\
-z+a_4+a_1\cdot a_5 (1-v)(1-w)\\
-v+a_6x+a_7\\
-w+(a_8-1)\left(1-\dfrac{z+a_9}{y+a_9} \right)\left(a_8+ \dfrac{z+a_9}{y+a_9}\right)^{-1}
\end{pmatrix},
\end{eqnarray}

then, solving the system \eqref{eq:004} is equivalent to finding a value $ X_\xi $ for  the function \eqref{eq:005} such that $\norm{f(X_\xi)}=0$.

\subsection{Generating an Initial Condition}

Without loss of generality, we may suppose that we have a function $ f: \Omega \subset \nset{R} \to \nset{R} $ with a simple root $ \xi $, that is,

\begin{eqnarray*}
f(x)=(x-\xi)g(x), & g(\xi)\neq 0,
\end{eqnarray*}

it should be noted that  in $ \Omega $ it is possible to find pairs of points $x_a$ and $x_b$, with $ x_a \neq x_b $, such that

\begin{eqnarray*}
f(x_a)\cdot f(x_b)\leq 0,
\end{eqnarray*}

as a consequence

\begin{eqnarray*}
\begin{array}{cccccc}
x_\xi \in [x_a,x_b]& \mbox{ or }  & x_\xi \in [x_b,x_a]& \mbox{ with }& f(x_\xi)=0.
\end{array}
\end{eqnarray*}

Hence, one way to approach the value $ x_\xi $, is to generate a set of $ N $ pseudorandom numbers $ \set {x_i}_{i = 1}^N $, with $ x_i <x_j \ \forall i <j $ and $ x_i \in \Omega \ \forall i \geq 1 $, with the intention of forming intervals $ [x_i, x_j] $ to evaluate the function $ f $ at its ends until finding one interval where it holds that

\begin{eqnarray*}
f(x_i)\cdot f(x_j)\leq 0,
\end{eqnarray*}

then, it is possible to take an initial condition $ x_0 \in [x_i, x_j] $ for use the iterative method \eqref{eq:c2.40}.

\subsection{Inspecting the System of Equations}

Assuming that the system \eqref{eq:004} has a solution, then it is possible to find two values $X_a$ and $X_b$, with $[X_a]_k\neq [X_b]_k \ \forall k\geq 1$, such that

\begin{eqnarray}\label{eq:006}
[f]_k(X_a)\cdot [f]_k(X_b)\leq 0, & \forall k\geq 1,
\end{eqnarray}

in consequence

\begin{eqnarray}\label{eq:007}
\begin{array}{cccc}
[X_\xi]_k\in \left[ [X_a]_k, [X_b]_k \right] & \mbox{or}&  [X_\xi]_k\in \left[ [X_b]_k, [X_a]_k \right] , & \forall k\geq 1.
\end{array}
\end{eqnarray}

To determine if the system \eqref{eq:004} has a solution, we may take the following values

\begin{eqnarray*}
X_a=(53,51,22,0,0)^T & \Rightarrow & f(X_a)\approx (1.831, 23.041,1.686,0.424,0.016)^T,\\
X_b=(54,52,23,1,1)^T & \Rightarrow & f(X_b)\approx (-2,-29,-3,-0.576,-0.984)^T,
\end{eqnarray*}

because the condition \eqref{eq:006} is satisfied, it is possible to guarantee that there is a value $ X_\xi $ that satisfies the condition \eqref{eq:007} with $X_\xi$ a zero of \eqref{eq:005}.

Although we have determined that the system \eqref{eq:004} has a solution, we need have in mind that, since the system  is nonlinear iterative methods as \eqref{eq:c2.40} are needed to try to solve it. Using iterative methods does not guarantee that the value $ X_\xi $ may be found. However, it is possible to find approximate values $X_N$ given as follows

\begin{eqnarray}\label{eq:009}
X_N=X_\xi+ \delta_\xi, & \delta_\xi=\delta_\xi( N), &\norm{\delta_\xi}<1.
\end{eqnarray}

Considering \eqref{eq:009}, it is necessary to give a general idea to determine if the function \eqref{eq:005}  is \textbf{stable} with respect to the  values $ X_N $.

\begin{definition}
Let $f:\Omega \subset \nset{R}^n \to \nset{R}^n$. If $ X_\xi $ is a zero of the function $ f $, we say that the function is stable with respect to the value $ X_\xi $, if when doing $ X_ \xi \to X_\xi +  \delta_\xi $, with $ \norm{\delta_\xi} <1 $, it holds that

\begin{eqnarray}\label{eq:010}
\norm{f(X_\xi+ \delta_\xi)}= \norm{\delta_f}<1.
\end{eqnarray}

\end{definition}

The condition \eqref {eq:010}, implies that for a function $ f $ to be stable, it is necessary that a slight perturbation $ \delta_ \xi $ in its solutions does not generate a great perturbation $ \delta_f $ in its images. To try to analyze the stability of the function \eqref{eq:005}, we may consider the following values

\begin{eqnarray*}
X_{N_1}=(53.8,51.6,22.1,0.3,0.1)^T & \Rightarrow & \norm{f(X_{N_1})}\approx 3.332, \\
X_{N_2}=(53.8,51.6,22.1,0.4,0.1)^T & \Rightarrow & \norm{f(X_{N_2})}\approx 1.408, \\
X_{N_3}=(53.8,51.6,22.1,0.5,0.1)^T & \Rightarrow & \norm{f(X_{N_3})}\approx 6.105, \\
\end{eqnarray*}

although $ X_{N_2} $ is not a zero of \eqref{eq:005}, it may be observed that it is a value close to the solution $ X_\xi $. Taking the canonical vector $ \hat{e} _4 = (0,0,0,1,0)^T $ and subtracting $ 1.408 $ in all the expressions on the right side, we obtain that

\begin{eqnarray*}
X_{N_1}=X_{N_2}-0.1\cdot\hat{e}_4 & \Rightarrow & \norm{f(X_{N_1})}-1.408\approx 1.924, \\
X_{N_2}=X_{N_2}+0.0\cdot\hat{e}_4 & \Rightarrow & \norm{f(X_{N_2})}-1.408\approx 0, \\
X_{N_3}=X_{N_2}+0.1\cdot\hat{e}_4 & \Rightarrow & \norm{f(X_{N_3})}-1.408\approx 4.697, \\
\end{eqnarray*}

the above helps us to visualize that a small perturbation $ \delta_{\xi} $ near the solution $ X_\xi $ is producing a large perturbation $ \delta_f $ near its image, reason why we can say that at first instance the function \eqref{eq:005} is unstable with respect to the values $X_{N_k}$. The mentioned above may be visualized in the Figure \ref{fig:01}.

\begin{figure}[!ht]
\begin{subfigure}{.5\textwidth}
\includegraphics[width=\textwidth, height=0.62\textwidth]{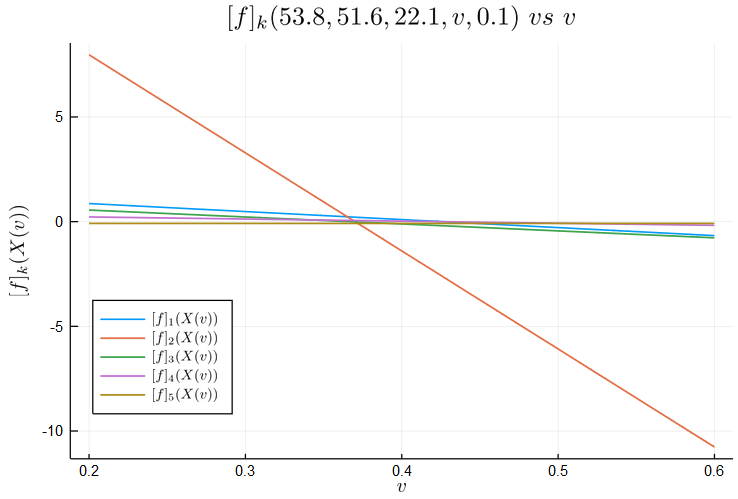}
\end{subfigure}
\begin{subfigure}{.5\textwidth}
\includegraphics[width=\textwidth, height=0.62\textwidth]{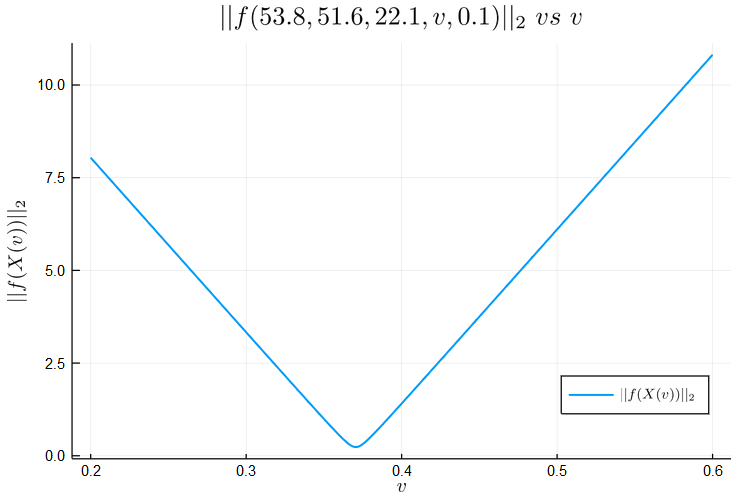}
\end{subfigure}
\caption{Graphs of the components $ [f]_k (X (v)) $ and $ \norm{f (X (v))}_2 $ with respect to different values of $ v $. }\label{fig:01}
\end{figure}

\subsection{Finding a Solution}

Considering that the function \eqref{eq:005} is unstable, it is necessary to give an initial condition $ X_ 0 $ very close to the solution $ X_\xi $ to be able to successfully use the iterative method  \eqref{eq:c2.40},  taking $X_0=X_{N_2}$ we obtain that

\begin{eqnarray}
f(X_0)\approx (0.098,-1.398,-0.11,0.024,-0.084)^T,
\end{eqnarray}

then, taking the fractional derivative given by \eqref{eq:c1.13}, the fractional pseudo-Newton method was implemented in the function \eqref{eq:005} to generate a sequence $\set{X_i}_{i=0}^N$ that approaches the solution $X_\xi$. In consequence, the results shown in the Table \ref{tab:01} were obtained.

\begin{table}[!ht]
\centering
$
\begin{array}{c}
\begin{array}{c|cccccc}
\toprule
&\alpha&x&y&z&v&w \\ \midrule
1&1.02934&53.80159759&51.59708283&22.09436105&  0.4243031& 0.01524 \\ \bottomrule
\end{array}
\\ \\
\begin{array}{ccc}
\toprule
||X_{N}-X_{N-1}||_2&||f(X_N)||_2&N \\ \midrule
2.04578e-3
&4.98732e-3
&  606\\ \bottomrule
\end{array}
\end{array}
$
\caption{Results obtained using the iterative method \eqref{eq:c2.40} with $\epsilon=e-4$. }\label{tab:01}
\end{table}

Using the values of the Table \ref{tab:01}, we obtain that

\begin{eqnarray}
f(X_N)\approx (0.001,0.0,-0.005,-0.0,0.001)^T,
\end{eqnarray}

the previous result allows us to get an idea of how close the solution $ X_N $, obtained by the method \eqref{eq:c2.40}, is to the solution $ X_\xi $ of the system \eqref{eq:004}.

\section{Conclusions}

The fractional pseudo-Newton method  to solve the problem of the need to invert a matrix in each iteration that is present in other methods. However, this method, may has at most an order of convergence (at least) linear,  and hence,  a speed of convergence relatively slow. As a consequence, it is necessary to use a larger number of positive values $ \alpha\in[0,2]\setminus \nset{Z}$  and a greater number of iterations, then we require a longer runtime to find solutions, so it may be considered as a method slow and costly,  but it is easy to implement.

This fractional iterative method may solve some nonlinear systems and linear systems and is efficient to find multiple solutions, both real and complex, using real initial conditions. It should be mentioned that this method is extremely recommended in systems that have infinite solutions or a large number of them.

\bibliography{Biblio}
\bibliographystyle{unsrt}

\nocite{rezania2016coupled}
\nocite{munoz2015efficiencies}

\end{document}